\documentclass[11pt]{article}
\usepackage{amsfonts,amssymb,eucal,a4}

\newtheorem{Abschnitt}{\vspace*{0pt}\hspace*{-0.8ex}}[section]

{\begin{Abschnitt}\label{#2}{\hspace{-1.1ex}\bf.~#1}%
}{\end{Abschnitt}}

\newenvironment{A}[2]%
{\begin{Abschnitt}\label{#2}{\hspace{-1.1ex}\bf.~#1}\hspace*{0ex}%
\rm%
}{\end{Abschnitt}}

\newcommand{\bA}{\begin{A}{}}
\newcommand{\eA}{\end{A}}

\newcommand{\bAZ}[2]{\begin{A}{#2.}{#1}}

\newcommand{\Formula}[2]{ \\[0.2cm] \hspace*{\fill} #1
                                    \hspace*{\fill} {\rm #2}\\[0.2cm]}

\newcommand{\proof}{\noindent \sc Proof.  \rm }%\small

\newcommand{\proofend}{\quad Q.\,E.\,D.}

\newcommand{\Z}{{\mathbb Z}}

\newcommand{\Zn}{{\Z}_n}

\renewcommand{\cal}[1]{{\EuScript #1}}

\begin{document}

%\date{April 10, 2001}

\title{\bf\protect Some Identities for Enumerators of Circulant Graphs}

\author{\protect{\bf Valery Liskovets}\thanks{
 Institute of Mathematics, National Academy of Sciences,
 220072, Minsk, Belarus,
 {\tt liskov{@}im.bas-net.by}
}}

\maketitle

\begin{abstract}
We establish analytically several new identities connecting enume\-ra\-tors
of different types of circulant graphs of prime, twice prime and
prime-squared orders. In particular, it is shown that
the semi-sum of the number of undirected circulants and the number of
undirected self-complementary circulants of prime order is equal to
the number of directed self-complementary circulants of the same order.

\medskip\noindent{\rm\small
}


\end{abstract}

\section{Introduction\label{Intro}}

Identities considered in this paper connect different
enumerators of circulant graphs mainly
of prime, twice prime and prime-squared orders. The idea of this paper
goes back to the article~\cite{KL96}, where we
counted uniformly circulants of five kinds and derived
several identities. Here we consider six types of circulants:
directed, undirected and oriented circulants (specified by valency
or not), and self-complementary circulants of the same types. Most of
the obtained identities may be called {\sl analytical} (or formal)
in the sense that they rest exclusively on the enumerative formulae
and follow from special properties of the cycle indices of regular
cyclic groups.
As a rule, it is more difficult to discover such an identity than to
prove it analytically. Almost all identities were first revealed and
conjectured due to numerical observations and computational experiments.

From the combinatorial point of view, most of the identities look rather
strange. They are very simple but no structural or algebraic properties
of circulants are used to derive them (with few exceptions), nor
establish we bijective proofs. The latter is challenging although
in some cases it is doubtful that there exist natural bijections
between participating circulants. Of course there may exist
other combinatorial or algebraic explanations or interpretations
of the identities.

Several identities hold only for a special type of prime orders $p$,
namely, those for which $\frac{p+1}{2}$ is also prime. Such primes are
familiar in number theory. Probably this is the first combinatorial
context where they play a substantial role.

We comprise here numerous identities that have been obtained previously
and deduce about ten new ones. We deliberately represent new identities
in different equivalent forms and formulate simple corollaries keeping
in mind possible future generalizations and combinatorial proofs.
Some derived identities look more elegant than the original ones.

The present paper is partially based upon the work~\cite{KL01} that contains
most detailed formulae for circulants, vast tables and several identities.
Here we reproduce all necessary results from it, and our exposition
is basically self-contained.

\bAZ{Int1}{Main definitions}
Let $n$ be a positive integer, $\Zn:=\{0,1,2,\dots,n-1\}$.
We denote by $\Zn^*$ the set of numbers in
$\Zn$ relatively prime to $n$ (that is invertible elements
{\sl modulo} $n$). So, $|\Zn^*|=\phi(n)$, where $\phi(n)$ is
the Euler totient function.
$Z(n)$ denotes a {\sl regular cyclic} permutation group of order
and degree $n$, i.e. the group gene\-rated by an $n$-cycle.

The {\it cycle index} of $Z(n)$ is the polynomial
\Formula{$\displaystyle{{\cal I}_n({\bf x})
 =\frac{1}{n}\sum\limits_{r|n}\phi(r)x_r^{n/r},}$}{(\ref{Intro}.1)}
where ${\bf x}$ stands for the sequence of variables $x_1,x_2,x_3,\dots$

The term ``graphs'' means both undirected and directed graphs. We
consider only simple graphs, that is graphs without loops and multiple
edges or arcs. %An $n$-graph means a graph of order $n$, where
The {\it order} of a graph means the number of its vertices.
We refer to Harary~\cite{Ha69} for notions concerning graphs.

An (undirected) edge is identified with the pair of the corresponding
oppositely directed arcs. Accordingly, an undirected graph is considered
as a (symmetric) digraph. On the contrary, a digraph is
{\it oriented} if it has no pair of oppositely directed arcs.

A circulant graph of order $n$, or simply a {\it circulant}, means a graph
$\Gamma$ on the vertex set $\Zn$ which is invariant with respect to the
cyclic permutation $(0,1,2,\dots,n-1)$, i.e.
if $(u,v)$ is an edge of $\Gamma$ then such is $(u+1,v+1)$.
In other terms, this is a Cayley graph with respect to the cyclic group
$\Zn$.
Every circulant is a regular graph of some valency $r$.

Graphs are considered here up to isomorphism. We deal with
the enumerators of (non-isomorphic) circulants of several types.
For convenience, the type is written as the subscript. Henceforth:
\begin{itemize}
 \item $C_{\rm d}(n)$ denotes the number of
{\bf d}irected circulant graphs;
 \item $C_{\rm u}(n)$ denotes the number of
{\bf u}ndirected circulant graphs;
 \item $C_{\rm o}(n)$ denotes the number of
{\bf o}riented circulant graphs;
 \item $C_{\rm sd}(n)$ and $C_{\rm su}(n)$ denote the numbers of
{\bf s}elf-complementary {\bf d}irected and
{\bf u}ndirected circulant graphs respectively;
 \item $C_{\rm t}(n)$ denotes the number of circulant {\bf t}ournaments;
  \item $C_{\rm d}(n,r)$, $C_{\rm u}(n,r)$ and $C_{\rm o}(n,r)$ denote
the corresponding numbers of circulants of order $n$ and valency $r$ while
$c_{\rm d}(n,z)$, $c_{\rm u}(n,z)$ and $c_{\rm o}(n,r)$ are their
generating functions by valency (polynomials in\ $z$):
\[\hspace*{-3em} c_{\rm d}(n,z):=\sum\limits_{r\geq 0}{C_{\rm d}(n,r)z^r},\quad
 c_{\rm u}(n,z):=\sum\limits_{r\geq 0}{C_{\rm u}(n,r)z^r},\quad
 c_{\rm o}(n,z):=\sum\limits_{r\geq 0}{C_{\rm o}(n,r)z^r}.\]
\end{itemize}
\indent
Clearly
$C_{\rm d}(n)=c_{\rm d}(n,z)|_{z:=1}=c_{\rm d}(n,1),\ C_{\rm u}(n)=c_{\rm u}(n,1)$
and $C_{\rm o}(n)=c_{\rm o}(n,1)$.

In more detail these quantities and the corresponding circulants
are consi\-de\-red in~\cite{KL96,KL01}. In particular,
the following simple uniform enumerative formulae have been obtained there:
\eA
%----------------------------
\bAZ{Int2}{Theorem (counting circulants of prime and twice prime order)}
{\it For any odd prime $p$,}
\begin{eqnarray*}
%1
c_{\rm d}(p,z)&=&{\cal I}_{p-1}({\bf x})|_{\{x_r:=1+z^r\ \}_{r=1,2,\ldots}}\\
%2
c_{\rm u}(p,z)&=&{\cal I}_{\frac{p-1}2}({\bf x})|_{\{x_r:=1+z^{2r}\}_{r=1,2,\ldots}}\\
%6
c_{\rm o}(p,z)&=&{\cal I}_{p-1}({\bf x})|_{\{x_r:=1\}_{r\ {\rm even}},\
    \{{x^2_r:=1+2z^r\}_{r\ {\rm odd}} }}\\
%4
C_{\rm sd}(p) &=&{\cal I}_{p-1}({\bf x})
 |_{\{{x_r:=2\}_{r\ {\rm even}},\ \{x_r:=0\}_{r\ {\rm odd}} }}\\
%5
C_{\rm su}(p) &=&{\cal I}_{\frac{p-1}2}({\bf x})
 |_{\{{x_r:=2\}_{r\ {\rm even}},\ \{x_r:=0\}_{r\ {\rm odd}} }}\\%%}}
%3
C_{\rm t}(p)  &=&{\cal I}_{p-1}({\bf x})|_{\{x_r:=0\}_{r\ {\rm even}},\
    \{{x^2_r:=2\}_{r\ {\rm odd}} }}\\
%7
c_{\rm d}(2p,z)&=&{\cal I}_{p-1}({\bf x})
 |_{\{x_r:=(1+z^r)^2\ \}_{r=1,2,\ldots}}\cdot(1+z)\\
%8
c_{\rm u}(2p,z)&=&{\cal I}_{\frac{p-1}2}({\bf x})
 |_{\{x_r:=(1+z^{2r})^2\}_{r=1,2,\ldots}}\cdot(1+z)\\
%9
c_{\rm o}(2p,z)&=&{\cal I}_{p-1}({\bf x})
 |_{\{x_r:=1\}_{r\ {\rm even}},\ \{{x_r:=1+2z^r\}_{r\ {\rm odd}} }}.
\end{eqnarray*}

\eA
%-----------------------------
\section{Cycle indices of cyclic groups}\label{I_m}

There are several technical formulae connecting the cycle indices
${\cal I}_{\frac{p-1}{2}}$ and ${\cal I}_{p-1}$. They are interesting
per se and will be used in the proofs of subsequent identities. % in their own right

For any natural $m$, we set \[m:=2^km'\] where $m'$ is odd.

In the polynomial ${\cal I}_{2m}$ we first distinguish the terms
corresponding to the divisors $r$ with the highest possible power of 2,
i.e. $k+1$:
\[{\cal I}_{2m}({\bf x})=\frac{1}{2m}\sum\limits_{r|2m}\phi(r)x_r^{2m/r}=
\frac{1}{2m}\Bigl(\sum\limits_{r|m}\phi(r)x_r^{2m/r}+
\sum\limits_{r|m'}\phi(2^{k+1}r)x_{2^{k+1}r}^{m'/r}\Bigr).\]
After easy transformations taking into account that $\phi(2^{k+1}r)=
2^k\phi(r)$ for odd $r$ and $k\geq 0$ we obtain
%------------------
\bAZ{Li1}{Lemma}
\Formula{$\displaystyle{2\,{\cal I}_{2m}({\bf x})={\cal I}_m({\bf x}^2)
+{\cal I}_{m'}({\bf x}_{(k+1)})}$}{(\ref{I_m}.1)}
{\it where ${\bf x}^2:=x_1^2,x_2^2,x_3^2,\dots$ and
${\bf x}_{(k+1)}:=x_{2^{k+1}},x_{2\cdot2^{k+1}},x_{3\cdot2^{k+1}}\dots$}

\medskip
Now in ${\cal I}_m({\bf x})$ we partition the set of divisors with respect
to powers of 2:
\[{\cal I}_m({\bf x})=
 \frac{1}{m}\Bigl(\sum\limits_{r|m'}\phi(r)x_r^{2^km'/r}+
 \sum\limits_{i=1}^k\sum\limits_{r|m'}2^{i-1}\phi(r)x_{2^ir}^{2^{k-i}m'/r}\Bigr)\]
and the same for ${\cal I}_{2m}({\bf x})$. Comparing similar terms in
both formulae, we easily arrive at the following:
\Formula{$\displaystyle{{\cal I}_m({\bf x})={\cal I}_{2m}(0,x_1,0,x_2,0,x_3,0,\dots)+
\frac{1}{2m}\sum\limits_{r|m'}\phi(r)x_r^{m/r}.}$}{(\ref{I_m}.2)}

The second summand on the right-hand side of formula~(\ref{I_m}.2) can
be represented in different useful forms. First of all, this is evidently
$\frac{1}{2m}\sum\limits_{r|m\atop r\ {\rm odd}}\phi(r)x_r^{m/r}.$
And this is also $\frac12{\cal I}_m(x_1,0,x_3,0,x_5,0,\dots)$. Hence
\Formula{$\displaystyle{{\cal I}_m({\bf x})={\cal I}_{2m}(0,x_1,0,x_2,0,x_3,0,\dots)+
\frac12{\cal I}_m(x_1,0,x_3,0,x_5,0,\dots).}$}{(\ref{I_m}.3)}
Every term in ${\cal I}_m$ contains only one variable. Therefore
\[{\cal I}_m(x_1,0,x_3,0,x_5,0,\dots)={\cal I}_m({\bf x})-
{\cal I}_m(0,x_2,0,x_4,0,x_6,0,\dots).\]
Hence by~(\ref{I_m}.3) we have
\eA
%------------------
\bAZ{Li2}{Lemma}
\Formula{$\displaystyle{2\,{\cal I}_{2m}(0,x_1,0,x_2,0,x_3,0,\dots)=
{\cal I}_m({\bf x})+{\cal I}_m(0,x_2,0,x_4,0,x_6,0,\dots)},$}{(\ref{I_m}.4)}
{\it that is},
\[2\,{\cal I}_{2m}({\bf y})|_{\{y_r:=0\}_{r\ {\rm odd}},\
 \{y_r:=x_{r/2}\}_{r\ {\rm even}}} = {\cal I}_m({\bf x}) +
 {\cal I}_m({\bf y}) |_{\{y_r:=0\}_{r\ {\rm odd}},\
 \{y_r:=x_r\}_{r\ {\rm even}}}.\]

\medskip
Now ${\cal I}_m(x_1,0,x_3,0,x_5,0,\dots)=
2\,{\cal I}_{2m}(\sqrt{x_1},0,\sqrt{x_3},0,\sqrt{x_5},0,\dots)$. Therefore
by~(\ref{I_m}.3),
\Formula{${\cal I}_m({\bf x})={\cal I}_{2m}(0,x_1,0,x_2,0,x_3,0,\dots)+
{\cal I}_{2m}(\sqrt{x_1},0,\sqrt{x_3},0,\sqrt{x_5},0,\dots).$}{(\ref{I_m}.5)}
Since the non-zero variables in both right-hand side summands alternate,
one may join them into a single cycle index. This gives rise to the
following expression:
\eA
%------------------
\bAZ{Li3}{Lemma}
\Formula{${\cal I}_m({\bf x})=
{\cal I}_{2m}(\sqrt{x_1},x_1,\sqrt{x_3},x_2,\sqrt{x_5},x_3,\dots).$}{(\ref{I_m}.6)}
{\it In other words,} $\ {\cal I}_m({\bf x})
={\cal I}_{2m}({\bf y})|_{\{y_r^2:={x_r}\}_{r\ {\rm odd}},\
 \{y_r:=x_{r/2}\}_{r\ {\rm even}}}$.
\eA
%-------------------

Finally we need one further formula. Substituting~(\ref{I_m}.5)
into~(\ref{I_m}.4) we obtain
\Formula{$\begin{array}{lll}
{\cal I}_{2m}(0,x_1,0,x_2,0,x_3,\dots)&=&
 {\cal I}_{2m}(\sqrt{x_1},0,\sqrt{x_3},0,\sqrt{x_5},\dots)\\
 &+&{\cal I}_m(0,x_2,0,x_4,0,x_6,\dots).
\end{array}$}{(\ref{I_m}.7)}

\vspace{-4ex}
%---------------------------------------
\section{Known identities}\label{Int}

\bA{Int3.1} Let $p$ be a prime such that $q=\frac{p+1}2$ is
also prime. Then by Klin -- Liskovets -- P\"oschel ~\cite{KL96},
\Formula{$c_{\rm u}(p,z)=c_{\rm d}(\frac{p+1}2,z^2),$}{(\ref{Int}.1)}
that is,
\Formula{$\quad C_{\rm u}(p,2r)=C_{\rm d}(\frac{p+1}2,r),\quad r\geq 0,$}{}
and
\Formula{$C_{\rm su}(p)=C_{\rm sd}(\frac{p+1}2).$}{(\ref{Int}.2)}
\indent
These equalities follow directly from Theorem~\ref{Int2} and are
in fact the first formal (i.e. analytically proved) identities
for enumerators of circulants.

We note that
\[p-1=2(q-1)\]
%partially
what explains a specific role of such primes in our considerations.

It follows from~(\ref{Int}.1) that
\Formula{$C_{\rm u}(p)=C_{\rm d}(\frac{p+1}2).$}{(\ref{Int}.$1'$)}
\vspace*{-4ex}
\eA
%--------------------------
\bA{Int3.2} If $p>3$ is a prime
such that $q=\frac{p+1}{2}$ is also prime, then
\Formula{$2c_{\rm o}(p,z) = c_{\rm o}(p+1,z)+1.$}{(\ref{Int}.3)}
\noindent{\sc Proof}~\cite{KL01}. Identity~(\ref{Int}.3) follows
directly from Theorem~\ref{Int2} (the third and ninth formulae) and
the equality
%Th1
\[2\,{\cal I}'_{2m}({\bf x}) = {\cal I}'_m({\bf x}^2)+1 \]
for an arbitrary $m$ where
${\cal I}'_m({\bf x}):={\cal I}_m({\bf x})|_{\{x_r:=1\}_{r\ {\rm even}}}$.
This equality is a parti\-cu\-lar case of expression~(\ref{I_m}.1)
since ${\cal I}_m(1,1,1,\dots)=1$.
Here we put $2m:=p-1$ (hence $m=q-1$).
\proofend

\medskip
Putting $z:=1$ we obtain
\Formula{$2C_{\rm o}(p) = C_{\rm o}(p+1)+1.$}{(\ref{Int}.$3'$)}
\vspace*{-4ex}
\eA
%----------------------------
\bA{Int3.3} According to~\cite{KL96},
\Formula{$C_{\rm su}(n)=0$}{(\ref{Int}.4)}
and
\Formula{$C_{\rm sd}(n)=C_{\rm t}(n)$}{(\ref{Int}.5)}
if $n=p$ or $p^2$ and $p\equiv 3$~(mod~4).

Next, combining~(\ref{Int}.5) with~(\ref{Int}.2) we obtain
\Formula{$C_{\rm su}(p)=C_{\rm t}(\frac{p+1}2)$}{(\ref{Int}.6)}
if both $p$ and $\frac{p+1}{2}$ are prime and $p\equiv 5$~(mod~8).

The least $p$ that meets the first two conditions and does not meet
the third is 73.
\eA
\bA{Int3.4} For any prime $p$,
\Formula{$\ C_{\rm sd}(p)=C_{\rm t}(p)+C_{\rm su}(p)$}{(\ref{Int}.7)}
\indent
Since tournaments and undirected self-complementary circulants
are particular cases of directed self-complementary circulants
(hence in general $\ C_{\rm sd}(n)\geq C_{\rm t}(n)+C_{\rm su}(n)$\,),
equality~(\ref{Int}.7) has a simple sense: any directed self-complementary
circulant graph of prime order is either anti-symmetric (a tournament) or
symmetric (an undirected graph).
This beautiful claim has been first established by Chia~--~Lim~\cite{CL86}
by means of simple algebraic arguments.
But in view of Theorem~\ref{Int2}
(the fourth, fifth and sixth formulae), identity~(\ref{Int}.7)
is a direct consequence of formula~(\ref{I_m}.7): merely substitute
2 for all variables $x_1,x_2,x_3,\dots$
\eA
%------------------
\bA{Int3.4a}
According to Fron\v{c}ek -- Rosa -- \v{S}i\'ra\v{n}~\cite{FR96} (see
also~\cite{AM99}), undirected self-com\-ple\-men\-tary circulants of
order $n$ exist if and only if all prime divisors $p$ of $n$ are
congruent to~1 modulo~4. Hence~(\ref{Int}.4) holds if there is a prime
$p|n,\ {p\equiv 3}$~(mod~4).
\eA
%-----------------
\bA{Int3.4b}
For composite orders, directed self-complementary circulants that are
neither tournaments nor undirected graphs do exist but are comparatively
rare. We call them {\it mixed}. The least suitable order is 15:\
$C^{\rm mixed}_{\rm sd}(15):=C_{\rm sd}(15)-C_{\rm su}(15)-C_{\rm t}(15)=
20-0-16=4$ (see Table~\ref{t1} in the Appendix). The four circulants are
constructed in~\cite{KL01}.
We will return to mixed circulants in Sections~\ref{M} and~\ref{Int3.4q}.
\eA
%-------------------------
\bA{Int3.5}
The last known identity concerns undirected circulants of even order
and odd valency:
\Formula{$C_{\rm u}(2n,2r+1)=C_{\rm u}(2n,2r)$}{(\ref{Int}.8)}
for any $n$ and $r$. This identity is known to hold for square-free $n$.
Moreover it has been verified for all orders less 54 and is conjectured
to be valid for all even orders (McKay~\cite{Mc95}; see also~\cite{KL01},
where a stronger conjecture concerning isomorphisms of circulants has been
proposed).
\eA
%---------------------------------

\section{New identities for circulants of prime order}\label{NId}

\bAZ{C3}{Proposition} {\it For prime $p$,
\Formula{$2C_{\rm sd}(p)=C_{\rm u}(p)+C_{\rm su}(p).$}{(\ref{NId}.1)}
In particular,
\Formula{$C_{\rm u}(p)=2C_{\rm sd}(p)
 =2C_{\rm t}(p)\quad {\rm if}\ \ p\equiv3\!\pmod 4.$}{(\ref{NId}.$1'$)}}
\proof
Substitute 2 for all variables in formula~(\ref{I_m}.4) with
$p-1=2m$. By Theorem~\ref{Int2} (the fourth, second and fifth formulae),
we immediately obtain~(\ref{NId}.1).
Clearly the second summand in~(\ref{I_m}.4) vanishes if $m$ is odd
(see~(\ref{Int}.4)).
%.
\proofend

\medskip
In Section~\ref{Ev} we will obtain a generalization of~(\ref{NId}.$1'$)
to $p\equiv1\!\pmod 4$.
\eA
%------------------
\bAZ{C3-R}{Remarks}
{\bf 1.} Despite that all participating quantities (and the corresponding
numerical values for small $p$) have been known long ago,
this striking identity has evidently escaped attention of the previous
researchers including the present author.

{\bf 2.} In view of equation~(\ref{Int}.7),
identity~(\ref{NId}.1) can be represented equivalently in the following form:
\Formula{$C_{\rm u}(p)=C_{\rm sd}(p)+C_{\rm t}(p)
 =C_{\rm su}(p)+2C_{\rm t}(p).$}{(\ref{NId}.$1''$)}
\indent
{\bf 3.} Is it possible to give a bijective proof of identity~(\ref{NId}.1)
or at least~(\ref{NId}.$1'$)? This question looks especially intriguing
in view of the fact that circulant graphs are naturally partitioned by
valency contrary to self-complementary circulants. Hence such a bijection
would introduce a certain artificial graduation (``pseudo-valency'') into
the class of self-complementary circulants of prime order. In particular,
some self-complementary graphs would correspond to the empty and complete
graphs.
To this end we could put formally $x_r:=1+z^r,\ r=1,2,\dots$,
in~(\ref{I_m}.4) instead of $x_r:=2$. But is there a natural combinatorial
interpretation of the coefficients of the left hand-side polynomial
thus obtained?

{\bf 4.} I do not know whether identity~(\ref{NId}.1) can be generalized
to non-prime orders.
\eA
%------------------------------
\bA{N11d2} %{Subtle identities}
We return to identity~(\ref{Int}.3).
There are subtler analogs of it for undirected and directed circulants.
By straightforward observations of numerical data and subsequent
numerical verifications with the help of the formulae for prime and
twice prime orders (Theorem~\ref{Int2}, the second, eighth, first and
seventh formulae) we arrived at the following somewhat unusual formulae:
\Formula{$4C_{\rm u}(p)
 =C_{\rm u}(p+1)+2\overline{C}_{\rm u}(2\tilde{p}+1),$}{(\ref{NId}.2)}
\vspace*{-4ex}
\Formula{$\displaystyle{2c_{\rm u}(p,z)=\frac{c_{\rm u}(p+1,z)}{1+z}
 +\overline{c}_{\rm u}(2\tilde{p}+1,z^{2^k})},$}{(\ref{NId}.3)}
\vspace{-3ex}
\Formula{$4C_{\rm d}(p)
 =C_{\rm d}(p+1)+2\overline{C}_{\rm u}(2\tilde{p}+1)$}{(\ref{NId}.4)}
and
\Formula{$\displaystyle{2c_{\rm d}(p,z)=\frac{c_{\rm d}(p+1,z)}{1+z}
 +\overline{c}_{\rm u}(2\tilde{p}+1,z^{2^{k}})}$}{(\ref{NId}.5)}
%%k+1?
whenever $p$ and $q=(p+1)/2$ are both odd primes.
Here $\tilde{p}$ denotes the maximal odd divisor of $p-1$~\footnote{In
the designations of Section~\ref{I_m}, $\tilde{p}=m'$ where
$p-1=2(q-1):=4m$.}
and $p-1:=2^{k+1}\tilde{p}$. Now
$\overline{c}_{\rm u}(2\tilde{p}+1,z):=c_{\rm u}(2\tilde{p}+1,z)$
if $2\tilde{p}+1$ is prime, otherwise $\overline{c}_{\rm u}$ is
calculated by {\sl the same} formula (the second formula in
Theorem~\ref{Int2}) despite that this time it
does not represent the number of non-isomorphic undirected circulants
of order $2\tilde{p}+1$.

\medskip
\proof
It is clear that formulae~{(\ref{NId}.2)} and~{(\ref{NId}.4)}
follow directly from~{(\ref{NId}.3)} and~{(\ref{NId}.5)} respectively.
The latter formulae are immediate consequences of equation~(\ref{I_m}.1)
with $q-1=2m$ and the corresponding formulae of Theorem~\ref{Int2}
for the orders $p$ and $p+1=2q$.
\proofend
\medskip

For instance, by data in Table~\ref{t2} one can verify that
$2c_{\rm d}(37,z)={c_{\rm d}(38,z)/(1+z)+c_{\rm u}(19,z^2)}$.
Hence for the valency $r=4$ we have numerically $2(1641+199)=3679+1$,
etc.

In particular, by~(\ref{NId}.3),
\Formula{$2C_{\rm u}(p,4r+2)=C_{\rm u}(p+1,4r+2)$}{(\ref{NId}.$3'$)}
since other terms correspond to undirected circulants of odd orders
and odd valency and, thus, vanish.

From~(\ref{NId}.2) and~(\ref{NId}.4) we obtain the following
identity not depending on $\tilde{p}$:
\Formula{$4C_{\rm d}(p)-C_{\rm d}(p+1)=
 4C_{\rm u}(p)-C_{\rm u}(p+1),\qquad \frac{p+1}2\ {\rm prime}.$}{(\ref{NId}.6)}
\indent
For example, for $p=13,\ 4\cdot352-1400=4\cdot14-48=8\ (=2C_{\rm u}(7))$.
For $p=73$ we obtain rather spectacularly
$4\cdot65588423374144427520-262353693496577709960=
4\cdot1908881900-7635527480=120\ (=2C_{\rm u}(19))$.\footnote{Moreover,
$120=4\cdot14602-58288=4C_{\rm u}(37)-C_{\rm u}(38)$.}

\medskip
Identity~(\ref{NId}.6) can also be written as
\Formula{$4(C_{\rm d}(p)-C_{\rm u}(p))=C_{\rm d}(p+1)-C_{\rm u}(p+1)$}{}
or
\Formula{$4C_{\rm d\backslash u}(p)
 =C_{\rm d\backslash u}(p+1),\qquad\frac{p+1}2\ {\rm prime},$}{(\ref{NId}.$6'$)}
where $C_{\rm d\backslash u}(n)$ denotes the number of directed
circulant graphs that are not undirected graphs.

Similarly from~(\ref{NId}.3) and~(\ref{NId}.5) we obtain
\Formula{$2(1+z)c_{\rm d\backslash u}(p,z)
 =c_{\rm d\backslash u}(p+1,z),\qquad\frac{p+1}2\ {\rm prime},$}{(\ref{NId}.7)}
or, equivalently,
\Formula{$2(C_{\rm d\backslash u}(p,r)+C_{\rm d\backslash u}(p,r-1))
 =C_{\rm d\backslash u}(p+1,r).$}{(\ref{NId}.$7'$)}
\indent
Thus, for example, for $p=13$ and $r=5$ we have
$C_{\rm d\backslash u}(13,5)=66-0=66,\ C_{\rm d\backslash u}(13,4)=43-3=40,
\ 66+40=106$ and $C_{\rm d\backslash u}(14,5)=217-5=2\cdot106$.
\eA
%---------------------------------
\bAZ{N11my}{Remark}
Can identities~(\ref{NId}.2) -- (\ref{NId}.7) (as well
as~(\ref{Int}.1) -- (\ref{Int}.3)) be treated bijectively? What then is
a sense of the sum or the corresponding difference? This is particularly
curious for~(\ref{NId}.2) and~(\ref{NId}.4) in the case of small
$\tilde{p}$. The existence of such a treatment seems doubtful
at least for composite $2\tilde{p}+1$. In this respect, identities
(\ref{NId}.$6'$) and~(\ref{NId}.$7'$) appear to be more promising.
\eA
%----------------------
%\newpage  %---------------- ad hoc ---------------
\bAZ{NT}{Number theoretic digression}
Some number theoretic aspects of identities~(\ref{NId}.2)~--~(\ref{NId}.7)
together with~(\ref{Int}.1)~--~(\ref{Int}.3) are worth considering.
There are 21 such pairs of primes $p=2q-1$ less 1000.
%and 14 between 1000 and 2000.
The first six $p$ are $3,5,13,37,61$ and 73 with their
corresponding $q =\ 2,3,7,19,31$ and 37. These are the sequences~M2492
and~M0849 in Sloane's Encyclopedia~\cite{SP95} (resp., A005383 and A005382
in its extended on-line version~\cite{SlXX}). In number theory these
numbers are called {\it nearly doubled primes}, and the pairs $q,p$ are
also known as {\it Cunningham chains of the second kind} of length~2
(see, e.g.,~\cite{Lo89,Fo99}). By definition, such primes $q$ resemble
the familiar Sophie Germain primes, that is, the primes $q$ such that
$p=2q+1$ is also prime. The latter primes play a different role in our
formulae: the polynomial ${\cal I}_{p-1}={\cal I}_{2q}$ contains
the minimal possible (for $p>3$) number of terms, four.

It is commonly believed that the set of nearly doubled primes
is infinite. Moreover, there is a conjecture that the number of
such primes $p<N$ grows asymptotically as $\frac{CN}{(\log{N})^2}$
where $C\dot=1.320$ (for $N=10^m$, this function is very close to
$\frac{10^m}{4m^2}$). Recall that the number of all primes $p<N$
grows approximately as $\frac{N}{\log{N}-1}$.

At present a lot of efforts in computational number theory are devoted
to the search for Cunningham chains of huge numbers, especially
long chains (see, e.g.,~\cite{Fo99}). In particular, the familiar program
proth.exe by Y.\,Gallot allows to effectively verify the primality of
numbers $\tilde{p}\cdot2^k+1$ with a fixed $\tilde{p}$. Keeping in mind
(\ref{NId}.2)~--~(\ref{NId}.5) we are especially interested in nearly
doubled primes with small $\tilde{p}$. In general it is easy to see
that such a pair $q,p$ can exist only if $3|\tilde{p}$. Here are the
current numerical results for $\tilde{p}\leq 27$.

%%$\tilde{p}=3$
Pairs of primes $q,p$ of the form $3\cdot2^k+1$
occur twice for $k\leq 303000$: only with $k=1,2$ and $k=5,6\ (p=193)$;
see the sequence~M1318 in~\cite{SP95} (or A002253~\cite{SlXX}).

%%$\tilde{p}=9$
Pairs of primes $q,p$ of the form $9\cdot2^k+1$
occur four times for $k\leq 145000$: with $k=1,2,\ k=2,3,\ k=6,7$
and $k=42,43$; see~M0751 (A002256).

%%$\tilde{p}=15$
Pairs of primes $q,p$ of the form $15\cdot2^k+1$
occur three times for $k\leq 184000$: with $k=1,2,\ k=9,10$ and $k=37,38$;
see~M1165 (A002258).

%%$\tilde{p}=21$
Pairs of primes $q,p$ of the form $21\cdot2^k+1$
occurs three times for $k\leq 164000$:
with $k=4,5,\ k=16,17$ and $k=128,129$ (see A032360~\cite{SlXX}).

%%$\tilde{p}=27$
Pairs of primes $q,p$ of the form $27\cdot2^k+1$
occurs twice for $k\leq 117000$:
with $k=19,20$ and $k=46,47$ (see A032363~\cite{SlXX}).
This gives rise to the least possible composite value of $2\tilde{p}+1$, 55.
So, for the first time it arises for $p=2q-1=27\cdot2^{20}+1=28311553$.

Clearly $2\tilde{p}+1=q$ if 8 does not divide $p-1$.
For $p<2000$, $2\tilde{p}+1$ turns out to be composite only in three cases.
$q=229,p=457$ is the least one; here $\tilde{p}=57$ and $2\tilde{p}+1=115$.

By numerical data we also found out that no Cunningham chain exists for
$\tilde{p}=51$ at least for $k<{140000}$. And the same for the numbers
$\tilde{p}=87$ and $\tilde{p}=93$.

Finally, two distinguished examples\footnote
{They are taken from the corresponding tables maintained in the WWW by
W.\,Keller and N.\,S.\,A.\,Melo, see
{\tt http://www.prothsearch.net/riesel.html} (cf. also~\cite{Ba79}).}:

$141\cdot2^k+1$ are prime for $k=555,556$;

$975\cdot2^k+1$ are prime for $k=6406,6407.$

\eA
%-----------------------------
\section{Circulants of prime-squared order}\label{M}

\bAZ{M2}{Theorem~{\rm\cite{KL96,KL01}}}{\it
\begin{eqnarray*}
%1
c_{\rm d}(p^2,z)&=&\ {\cal C}(p^2;{\bf x,y})
 |_{\{x_r:=1+z^r,\ \ y_r:=1+z^{pr}\ \}_{r=1,2,\ldots}}\\
%2
c_{\rm u}(p^2,z)&=& {\cal C}^*(p^2;{\bf x,y})
 |_{\{x_r:=1+z^{2r},\ y_r:=1+z^{2pr}\}_{r=1,2,\ldots}}\\
%6
c_{\rm o}(p^2,z)&=&\ {\cal C}(p^2;{\bf x,y})
 |_{\{x_r:=1,\ y_r:=1\}_{r\ {\rm even}},\
 \{{x^2_r:=1+2z^r,\ y^2_r:=1+2z^{pr}\}_{r\ {\rm odd}} }}\\
%4
C_{\rm sd}(p^2) &=&\ {\cal C}(p^2;{\bf x,y})
 |_{\{{x_r:=2,\ y_r:=2\}_{r\ {\rm even}},\ \{x_r:=0,\ y_r:=0\}_{r\ {\rm odd}} }},\\
%5
C_{\rm su}(p^2) &=& {\cal C}^*(p^2;{\bf x,y})
 |_{\{{x_r:=2,\ y_r:=2\}_{r\ {\rm even}},\ \{x_r:=0,\ y_r:=0\}_{r\ {\rm odd}} }}\\
%3
C_{\rm t}(p^2)  &=&\ {\cal C}(p^2;{\bf x,y})
 |_{\{x_r:=0,\ y_r:=0\}_{r\ {\rm even}},\
 \{{x^2_r:=2,\ y^2_r:=2\}_{r\ {\rm odd}} }} %.
\end{eqnarray*}
where
\[{\cal C}(p^2;{\bf x,y}):=\frac{1}{p}{\cal I}_{p-1}({\bf x}^{p+1})
 -\frac{1}{p}{\cal I}_{p-1}({\bf xy})
 +{\cal I}_{p-1}({\bf x}){\cal I}_{p-1}({\bf y})\]
and
\[{\cal C}^*(p^2;{\bf x,y})
 :=\frac{1}{p}{\cal I}_{\frac{p-1}2}({\bf x}^{p+1})
 -\frac{1}{p}{\cal I}_{\frac{p-1}2}({\bf xy})
 +{\cal I}_{\frac{p-1}2}({\bf x}){\cal I}_{\frac{p-1}2}({\bf y})\]
with ${\bf x}^{p+1}:=x_1^{p+1},x_2^{p+1},x_3^{p+1},\dots$ and
${\bf xy}:=x_1y_1,x_2y_2,x_3y_3,\dots$}
\eA
%-----------------
\bAZ{Mix}{Mixed self-complementary circulant graphs}
By definition (see Section~\ref{Int3.4b}),
\Formula{$C^{\rm mixed}_{\rm sd}(p^2)
 :=C_{\rm sd}(p^2)-C_{\rm su}(p^2)-C_{\rm t}(p^2).$}{(\ref{M}.1)}
According to~\cite{LP00,KL01}, the number of non-CI (non-Cayley isomorphic)
circulants of order $p^2$ is
\Formula{$D_{\rm i}(p^2)=C_{\rm i}(p)^2,$}{(\ref{M}.2)}
where ${\rm i}\in\{\rm sd,su,t\}$. We recall that a circulant is said
to be {\it non-CI} if there exists a circulant isomorphic but not Cayley
isomorphic to it. A {\it Cayley isomorphism} means an isomorphism that
is induced by an automorphism of the underlying group $\Z_n$.
\eA
%-------------------------------
\bAZ{M3}{Proposition}
\Formula{$C^{\rm mixed}_{\rm sd}(p^2)
 =2C_{\rm su}(p)C_{\rm t}(p)$}{(\ref{M}.3)}
and
\Formula{$C^{\rm mixed}_{\rm sd}(p^2)
 =D_{\rm sd}(p^2)-D_{\rm su}(p^2)-D_{\rm t}(p^2),$}{(\ref{M}.4)}
{\it that is, the mixed self-complementary circulants of order $p^2$
are exactly the \mbox{\rm non-CI} mixed self-complementary circulants.}

\medskip
\proof
We make use of an algebraic property of self-complementary circulants
of prime-power order. According to a result announced by Li~\cite{Li98}
(Theorem~3.3), if $\Gamma$ is a self-complementary circulant of order
$p^2$ then one of the following holds.
\begin{itemize}
\item $\Gamma$ can be obtained by means of
the well-known (alternating cycle) construction
discovered by Sachs and Ringel.
\item $\Gamma=\Gamma_1[\Gamma_2]$ %%\circ
where $\Gamma_1$ and $\Gamma_1$ are self-complementary circulants
of order $p$. Here $\Gamma_1[\Gamma_2]$ is the composition (called
also the wreath or lexicographic product) defined as follows:
in $\Gamma_1$ we replace each vertex by a copy of $\Gamma_2$;
each edge of $\Gamma_1$ gives rise to the edges connecting all
pairs of vertices from the two corresponding copies of $\Gamma_2$.
\end{itemize}
The first construction generates only undirected circulants or
tournaments (cf.~\cite{LP00}); moreover, all of them are CI.
Now, there is no mixed self-complementary circulant of order $p$
(this is identity~(\ref{Int}.7)).
Therefore the se\-cond construction gives rise to a mixed graph
if and only if one of the factors is an undirected
self-complementary circulant and the other factor is a tournament.
This proves~(\ref{M}.3). Further, all self-complementary circulants
$\Gamma=\Gamma_1[\Gamma_2]$ are non-CI~\cite{LP00}. This, together
with~(\ref{M}.2), proves~(\ref{M}.4) (moreover, this proves~(\ref{M}.2)
since the composition of two undirected circulants
is undirected and the composition of two tournaments is a tournament).
\proofend

It would be interesting to find an analytical derivation of
these equations with the help of Theorem~\ref{M2}.

\medskip
By~(\ref{M}.1) we have
\eA
%------------------------
\bAZ{M3Cor}{Corollary}
\Formula{$C_{\rm sd}(p^2)-C_{\rm su}(p^2)-C_{\rm t}(p^2)
 =C_{\rm sd}(p)^2-C_{\rm su}(p)^2-C_{\rm t}(p)^2.$}{(\ref{M}.5)}
%3
\vspace*{-4ex}
\eA
%--------------
\bAZ{M3Ex}{\bf Example} $p=13$. By Theorem~\ref{M2},\\
$C_{\rm sd}(13^2) =123992391755402970674764$,\
$C_{\rm su}(13^2) =56385212104$ and\\
$C_{\rm t}(13^2)\ =123992391755346585462636$.\
It follows that $C^{\rm mixed}_{\rm sd}(13^2)=24.$
Now $C_{\rm sd}(13)^2=8^2=64,\ C_{\rm su}(13)^2=2^2=4,\
C_{\rm t}(13)^2=6^2=36$ and $64-4-36=24=2\cdot2\cdot6.$
\eA

By~(\ref{Int}.7) (or, instead, by~(\ref{M}.1) and~(\ref{M}.3)),
identity~(\ref{M}.5) can be represented as follows:
%--------------------------
\Formula{$C_{\rm sd}(p^2)
 =C_{\rm su}(p^2)+C_{\rm t}(p^2)+2C_{\rm su}(p)C_{\rm t}(p).$}{(\ref{M}.6)}
\indent
We note also that if ${p\equiv 3}$~(mod~4), then $C_{\rm su}(p)$
and $C_{\rm su}(p^2)$ vanish by~(\ref{Int}.4), and identity~(\ref{M}.6)
turns into~(\ref{Int}.5) for $n=p^2$.

%-----------------------
\section{Alternating sums}\label{Alt}

Alternating sums serve as one further source of formal identities.
First consider
directed circulants of prime order. Take the generating function
$c_{\rm d}(p,t)$ and put $t:=-1$.
By Theorem~\ref{Int2} we see that the result is equal to
$C_{\rm sd}(p)$. By Theorem~\ref{M2}, the same result is valid for
the orders $n=p^2$. Moreover, by formulae given in~\cite{KL01}
this is valid for arbitrary odd square-free orders.
Therefore we have
\Formula{$c_{\rm d}(n,-1)=C_{\rm sd}(n).$}{(\ref{Alt}.1)}
\indent
The corresponding result holds for undirected circulants with respect
to the substitution $t^2:=-1$, or $t:=\sqrt{-1}$:
\Formula{$c_{\rm u}(n,t)|_{t^2:=-1}=C_{\rm su}(n).$}{(\ref{Alt}.2)}
Both formulae have been proved for square-free and prime-squared $n$
and it is natural to suggest that they are valid in general:

\bAZ{Conj1}{Conjecture} Identities~(\ref{Alt}.1) and~(\ref{Alt}.2)
hold for any odd order $n$.

\medskip
Trivially (by complementarity), identity~(\ref{Alt}.1) holds also for
even $n$, and~{(\ref{Alt}.2)} holds for $n\equiv 3$~(mod~4).
Identity~{(\ref{Alt}.2)} is also valid for $n=45$ as numerical
data~\cite{Mc95} show.
\eA

The behaviour of oriented circulant graphs is different.
Numerical observations show that
\Formula{$c_{\rm o}(n,-1)=0$}{(\ref{Alt}.3a)}
if $n$ has at least one prime divisor $p\equiv 3$~(mod~4), otherwise
\Formula{$c_{\rm o}(n,-1)=1.$}{(\ref{Alt}.3b)}
These identities hold for prime $n=p$ by Theorem~\ref{Int2},
for {\sl odd} square-free $n$ by~\cite{KL01} and for $n=p^2$
by Theorem~\ref{M2}. Again we {\bf conjecture} they to be valid
for all odd $n$.

For {\sl even} square-free $n$ we found that identity~(\ref{Alt}.3b)
holds if
$n=2n',\ n'$ odd, and (\ref{Alt}.3a) holds if $n=4n',\ n'$ square-free.
The behaviour of $c_{\rm o}(n,-1)$ for $n=8n',\ n'>1$, remains
unknown.

Identities~(\ref{Alt}.1) and~(\ref{Alt}.2) for prime $n=p$
transform~{(\ref{NId}.1)} into the following equality:
\Formula{$2c_{\rm d}(p,-1)
 =c_{\rm u}(p,1)+c_{\rm u}(p,\sqrt{-1}).$}{(\ref{Alt}.4)}
\vspace*{-4ex}
\bAZ{Ev}{Even- and odd-valent circulants}
Due to~(\ref{Alt}.1) and~(\ref{Alt}.2) we can find simple
expressions for the numbers of circulants of (non-specified) {\sl even}
(and, resp., {\sl odd}) valency; for undirected circulants we consider
only odd orders and mean even and odd {\sl semi}-valencies,
that is, valencies congruent, respectively, to 0 and 2 modulo~4.
We use the superscript {\small e} and {\small o} to denote these
numbers. Now, formula~(\ref{Alt}.1) is nothing than
$C^{\rm e}_{\rm d}(n)-C^{\rm o}_{\rm d}(n)=C_{\rm sd}(n)$. Since
$C^{\rm e}_{\rm d}(n)+C^{\rm o}_{\rm d}(n)=C_{\rm  d}(n)$,
we obtain
\Formula{$\displaystyle{C^{\rm e}_{\rm d}(n)
 =\frac{C_{\rm d}(n)+C_{\rm sd}(n)}2}$}{(\ref{Alt}.5e)}
and
\Formula{$\displaystyle{C^{\rm o}_{\rm d}(n)
 =\frac{C_{\rm d}(n)-C_{\rm sd}(n)}2.}$}{(\ref{Alt}.5o)}
\indent
So, these expressions hold for square-free, prime-squared and even $n$
and are assumed to hold for all orders.

Similarly,~(\ref{Alt}.2) gives rise to
\Formula{$\displaystyle{C^{\rm e}_{\rm u}(n)
 =\frac{C_{\rm u}(n)+C_{\rm su}(n)}2,}$}{(\ref{Alt}.6e)}
and
\Formula{$\displaystyle{C^{\rm o}_{\rm u}(n)
 =\frac{C_{\rm u}(n)-C_{\rm su}(n)}2}$}{(\ref{Alt}.6o)}
for undirected circulants of odd orders and even and, resp., odd
semi-valency. Equations~(\ref{Alt}.6e) and~(\ref{Alt}.6o) remain
unproven unless $n$ is square-free or prime-squared or congruent to~3
modulo~4.

Clearly the respective expressions can be extracted from~(\ref{Alt}.3a)
and~(\ref{Alt}.3b) for oriented circulants.

Comparing formula~(\ref{Alt}.6e) for prime $n=p$ with~(\ref{NId}.1)
we obtain the following curious identity:
\Formula{$C^{\rm e}_{\rm u}(p)=C_{\rm sd}(p).$}{(\ref{Alt}.7)}
This equation directly generalizes identity~(\ref{NId}.1$'$) to
$p\equiv 1$~(mod~4) because $C^{\rm e}_{\rm u}(p)=C_{\rm u}(p)/2$ for
$p\equiv 3$~(mod~4).
\eA
%-------------------
\section{Conclusion}\label{Fin}

We conclude that the enumerative theory of circulants is full
of hidden inter-dependencies, part of which are presented in this paper.
Table~\ref{t3} in the Appendix contains a summary of previous and new
identities.

We expect that there should exist further generalizations of the obtained
identities for other classes of circulant graphs, first of all, for
multigraphs and graphs with coloured or marked edges.

\bA{Fin1} In general, analytical identities are characteristic for the
enumerators of self-complementary graphs of diverse classes and can be
found in numerous publications. These results are collected in the
surveys by Robinson~\cite{Ro81} and Farrugia~\cite{Fa99} (mainly in Ch.\,7).
In the latter paper, several open questions are also posed. In particular,
the problem~{\sf K} in Sect.\,7.64 is just the problem of finding a natural
bijection for identity~(\ref{Int}.2).
\eA
%--------------------
\bAZ{Int3.4q}{Open question}
Is identity~(\ref{Int}.5) valid for the orders all whose prime divisors
are congruent to~3 modulo~4? In other words (since~(\ref{Int}.4) holds
according to Section~\ref{Int3.4a}), are there mixed self-complementary
circulants of such orders?
As conjectured in~\cite{KL01}, mixed self-complementary circulants of
order $n$ exist if and only if $n$ is odd composite and has a prime divisor
$p\equiv 1$~(mod~4). If so, then identity~(\ref{Int}.7) holds exactly for
the other orders. This claim is valid for square-free orders, and it can
also be proved for the prime-power orders $n=p^k$. %%was?
\eA
%------------------
\bA{Fin2} Identities~(\ref{Alt}.1),~(\ref{Alt}.2), (\ref{Alt}.5)
and~(\ref{Alt}.6) are rather typical; cf., e.g., my paper~\cite{Lis00},
where other examples of even- and odd-specified quantities and the
corresponding semi-sum expressions for them are given.
\eA
%-----------------
\bA{Fin3} Finally, instead of equalities, we touch one important type
of inequa\-li\-ties which are usually proved analytically. I suppose that
the sequence of the numbers $C_{\rm u}(p,2r),\ 1<r<(p-1)/2$, is
{\it logarithmically concave}, that is
\[C_{\rm u}(n,2r)^2\geq C_{\rm u}(n,2r-2)C_{\rm u}(n,2r+2)\]
for any prime order $n=p$ and $1<r<(n-1)/2$. In other words, the
sequence of ratios $C_{\rm u}(p,2r)/C_{\rm u}(p,2r+2)$ is increasing
except for the first and the last member. For composite orders this
does not necessarily hold.
In particular, the opposite inequality holds for $r=2$ when $n=27,\ 121$
and $169$. However I do not know counterexamples for square-free orders.
\eA

\bibstyle{plain}

\newpage

%------------------------------------------------
{
%%\small\tabcolsep=1.1ex
\begin{table}[htb]
\centerline{\large\bf Appendix: numerical results and summary}\label{A}

\medskip
{\normalsize Tables~\ref{t1} and~\ref{t2} contain relevant numerical data
obtained by Theorems~\ref{Int2} and~\ref{M2} (they partially reproduce
data from~\cite{KL01}).}

\caption{Non-isomorphic circulant graphs\label{t1}}
\medskip\small
\begin{tabular}{|r||r|r|r||r|r|r|}
\hline
$n$&$C_{\rm d}(n)$&$C_{\rm u}(n)$&$C_{\rm o}(n)$&$C_{\rm sd}(n)$
   &$C_{\rm su}(n)$&$C_{\rm t}(n)$\\
\hline
%1&       1      &     1 &1          &   1  &  1&    1  \\%
 2&       2      &     2 &1          &   0  &  0&    0  \\%
 3&       3      &     2 &2          &   1  &  0&    1  \\%
 4&       6      &     4 &2          &   0  &  0&    0  \\%
 5&       6      &     3 &3          &   2  &  1&    1  \\%
 6&      20      &     8 &5          &   0  &  0&    0  \\%
 7&      14      &     4 &6          &   2  &  0&    2  \\%
 8&      46      &    12 &7          &   0  &  0&    0  \\%McK
 9&      51      &     8 &16         &   3  &  0&    3  \\%
10&     140      &    20 &21         &   0  &  0&    0  \\%
11&     108      &     8 &26         &   4  &  0&    4  \\%
12&     624      &    48 &64         &   0  &  0&    0  \\%
13&     352      &    14 &63         &   8  &  2&    6  \\%
14&    1400      &    48 &125        &   0  &  0&    0  \\%
15&    2172      &    44 &276        &  20  &  0&   16  \\%
%16&             &    84 &           &   0  &  0&    0  \\%
17&    4116      &    36 &411        &  20  &  4&   16  \\%
18&   22040      &   192 &1105       &   0  &  0&    0  \\%
19&   14602      &    60 &1098       &  30  &  0&   30  \\%
20&   68016      &   336 &2472       &   0  &  0&    0  \\%
21&   88376      &   200 &4938       &  88  &  0&   88  \\%
22&  209936      &   416 &5909       &   0  &  0&    0  \\%
23&  190746      &   188 &8054       &  94  &  0&   94  \\%
%24&             &  1312 &           &   0  &  0&    0 \\%A=1344
25&  839094      &   423 &26577      & 214  &  7&  205  \\%
%\hline
26& 2797000      &  1400 &44301      &   0  &  0&    0  \\%
%27&             &   928 &           &      &  0&       \\%A_u=944
28&11276704      &  3104 &132964     &   0  &  0&    0  \\%
29& 9587580      &  1182 &170823     & 596  & 10&  586  \\%
30&67195520      &  8768 &597885     &   0  &  0&    0  \\%
31&35792568      &  2192 &478318     &1096  &  0& 1096  \\%1091
%32&             &  8364 &           &   0  &  0&    0  \\%A=8784
33&214863120     &  6768 &2152366    &3280  &  0& 3280  \\%
34&536879180     & 16460 &2690421    &   0  &  0&    0  \\%
35&715901096     & 11144 &5381028    &5560  &  0& 5472  \\%
%36&             & 46784 &           &   0  &  0&    0  \\%A=46848
37&1908881900    & 14602 &10761723   &7316  & 30& 7286  \\%7280
38&7635527480    & 58288 &21523445   &   0  &  0&    0  \\%
39&11454711464   & 44424 &48427776   &21944 &  0& 21856 \\%
%40&             &136128 &           &   0  &  0&    0  \\%A=138432
41&27487816992   & 52488 &87169619   &26272 & 56&26216  \\%26214
42&183264019200  &355200 &290566525  &   0  &  0&    0  \\%
43&104715443852  & 99880 &249056138  &49940 &  0& 49940 \\%49929
44&440020029120  &432576 &523020664  &   0  &  0&    0  \\%
%45&             &351424 &           &      &  0&       \\%+288
46&1599290021720 &762608 &1426411805 &   0  &  0&    0 \\%
47&1529755490574 &364724 &2046590846 &182362&  0&182362\\%182361
%48&             &2122944&           &   0  &  0&    0 \\%
49&6701785562464 &798952 &6724513104 &399472&  0&399472\\%
50&28147499352824&3356408&14121476937&   0  &  0&    0 \\%
%51&       .     &2105420&         . &     .&  0&    . \\%
%53&       .     &  .    &         . &     .&316&    . \\%
\hline
\end{tabular}
\end{table} }

{
%%\footnotesize
\tabcolsep=1.19ex
\begin{table}[htb]
\caption{Enumeration of circulant graphs by valency (for selective orders)\label{t2}} %for nearly doubled primes
\medskip\small
\begin{tabular}{|r||r|r|r|r|r|r|r|r|r|r|}
%\hline
\multicolumn{11}{c}{$c_{\rm u}(n,r),\ r$\ even}\\
\hline
&\multicolumn{10}{|c|}{$n$}\\
\hline
$r$&7&13&14&19&   37&  38&      61&        62&        73&        74\\
\hline
0 &1& 1& 1&  1&    1&   1&       1&         1&         1&         1\\
2 &1& 1& 2&  1&    1&   2&       1&         2&         1&         2\\
4 &1& 3& 5&  4&    9&  17&      15&        29&        18&        36\\
6 &1& 4& 8& 10&   46&  92&     136&       272&       199&       398\\
8 & & 3& 5& 14&  172& 340&     917&      1827&      1641&      3281\\
10& & 1& 2& 14&  476& 952&    4751&      9502&     10472&     20944\\
12& & 1& 1& 10& 1038&2066&   19811&     39591&     54132&    108264\\
14& &  &  &  4& 1768&3536&   67860&    135720&    231880&    463760\\
16& &  &  &  1& 2438&4862&  195143&    390195&    840652&   1681300\\
18& &  &  &  1& 2704&5408&  476913&    953826&   2615104&   5230208\\
20& &  &  &   & 2438&4862& 1001603&   2003005&   7060984&  14121968\\
22& &  &  &   & 1768&3536& 1820910&   3641820&  16689036&  33378072\\
24& &  &  &   & 1038&2066& 2883289&   5766243&  34769374&  69538738\\
26& &  &  &   &  476& 952& 3991995&   7983990&  64188600& 128377200\\
28& &  &  &   &  172& 340& 4847637&   9694845& 105453584& 210907168\\
30& &  &  &   &   46&  92& 5170604&  10341208& 154664004& 309328008\\
32& &  &  &   &    9&  17& 4847637&   9694845& 202997670& 405995326\\
34& &  &  &   &    1&   2& 3991995&   7983990& 238819350& 477638700\\
36& &  &  &   &    1&   1& 2883289&   5766243& 252088496& 504176992\\
38& &  &  &   &     &    & 1820910&   3641820& 238819350& 477638700\\
40& &  &  &   &     &    & 1001603&   2003005& 202997670& 405995326\\
%42& & &  &   &     &    &  476913&    953826& 154664004& 309328008\\
%44& & &  &   &     &    &  195143&    390195& 105453584& 210907168\\
%46& & &  &   &     &    &   67860&    135720&  64188600& 128377200\\
\hline
\end{tabular}

\tabcolsep=0.94ex
%\caption{Directed and oriented circulant graphs by valency (extraction)\label{t3}} %for nearly doubled primes
\medskip
\begin{tabular}{|r||r|r|r|r|r|r|r||r|r|r|r|}
%\hline
\multicolumn{8}{c||}{$c_{\rm d}(n,r)$}&\multicolumn{4}{c}{$c_{\rm o}(n,r)$}\\
\hline
&\multicolumn{7}{c||}{$n$}&\multicolumn{4}{c|}{$n$}\\
\hline
$r$&7&13& 14&  19&      31&        37&       38&13&14&     37& 38\\
\hline
0 &1& 1 &1  &1   &       1&1         &1        &1 &1 &1      &1\\
1 &1& 1 &3  &1   &       1&1         &3        &1 &2 &1      &2\\
2 &3& 6 &14 &9   &      15&18        &38       &5 &10&17     &34\\
3 &4& 19&50 &46  &     136&199       &434      &14&28&182    &364\\
4 &3& 43&123&172 &     917&1641      &3679     &20&40&1360   &2720\\
5 &1& 66&217&476 &    4751&10472     &24225    &16&32&7616   &15232\\
6 &1& 80&292&1038&   19811&54132     &129208   &6 &12&33006  &66012\\
7 & & 66&292&1768&   67860&231880    &572024   &  &  &113152 &226304\\
8 & & 43&217&2438&  195143&840652    &2145060  &  &  &311168 &622336\\
9 & & 19&123&2704&  476913&2615104   &6911508  &  &  &691494 &1382988\\
10& & 6 &50 &2438& 1001603&7060984   &19352176 &  &  &1244672&2489344\\
11& & 1 &14 &1768& 1820910&16689036  &47500040 &  &  &1810432&3620864\\
12& & 1 &3  &1038& 2883289&34769374  &102916810&  &  &2112184&4224368\\
13& &   &1  &476 & 3991995&64188600  &197915938&  &  &1949696&3899392\\
14& &   &   &172 & 4847637&105453584 &339284368&  &  &1392640&2785280\\
15& &   &   &46  & 5170604&154664004 &520235176&  &  &742752 &1485504\\
16& &   &   &9   & 4847637&202997670 &715323334&  &  &278528 &557056\\
17& &   &   &1   & 3991995&238819350 &883634026&  &  &65536  &131072\\
18& &   &   &1   & 2883289&252088496 &981815692&  &  &7286   &14572\\
19& &   &   &    & 1820910&238819350 &981815692&  &  &       &\\
20& &   &   &    & 1001603&202997670 &883634026&  &  &       &\\
%21& &  &   &    &  476913&154664004 &715323334&  &  &       &\\
%22& &  &   &    &        &105453584 &520235176&  &  &       &\\
%23& &  &   &    &        &64188600  &339284368&  &  &       &\\
\hline
\end{tabular}
\end{table}

{ %%\tabcolsep=1.0ex  %%\normalsize
\begin{table}[htb]
\caption{Systematized list of identities\label{t3}}
\medskip \small
\begin{tabular}{|r||c|c|c|c|l|c|}
\hline
No&\small Formula&Orders\,$\rm^a$&Restrictions&Types&\ \quad Proof&Refer. \\%
\hline
\multicolumn{7}{|c|}{For self-complementary circulants:}\\
\hline
 1&(3.4)  &$n$ &$\exists\,p|n,p\equiv3$~(4)&su &\quad Combin.&\cite{FR96}\\%5
  &       &       &                  &         &\quad Algebr.&\cite{AM99}\\%
  &       &$p$ or $p^2$&$p\equiv3$~(mod~4)\ &    &\quad Analyt.&\cite{KL96}\\%
 2&(3.5)  &$p$ or $p^2$&$\ p\equiv3$~(mod~4)\,$\rm^b$&t, sd&\quad Analyt.&\cite{KL96}\\%6
  &       &            &                            &     &$\Leftarrow$(3.4),\,(5.3)&(New)\\%
%%9&(3.7)  &$p$    &$p\equiv3$~(mod~4)&t, sd    &$\Leftarrow$(3.6),\,(3.4)&(New) \\%
 3&(3.7)  &$p$    &-                 &su, t, sd&\quad Algebr.&\cite{CL86}\\%6.1
  &       &$n$    &$n=p$,\dots $\rm^b$&        &\quad Analyt.&New        \\%
 4&(4.1)  &$p$    &-                 &u, su, sd&\quad Analyt.&New        \\%2
 5&(4.$1'$)&$p$   &$p\equiv3$~(mod~4)&u, sd    &$\Leftarrow$(4.1),\,(3.4)&(New) \\%3
 6&(4.$1''$)&$p$  &-                 &u, su, t &$\Leftarrow$(4.1),\,(3.7)&(New) \\%4
 7&(5.5)  &$p,p^2$&-                 &su, t, sd&$\Leftarrow$(5.2),\,(5.4)&(New) \\
 8&(5.6)  &$p,p^2$&-                 &su, t, sd&$\Leftarrow$(3.7),\,(5.5)&(New) \\
 9&(3.6)&$p,q$&$p\!+\!1\!=\!2q\equiv6$\,(8)&su, t&$\Leftarrow$(3.2),\,(3.7)&(New)\\%8
10&(3.2)  &$p,q$  &$p+1=2q$          &su, sd   &\quad Analyt.&\cite{KL96}\\%
\hline
\multicolumn{7}{|c|}{Other valency independent:}\\
\hline
11&(3.$1'$)&$p,q$ &$p+1=2q$          &u, d     &$\Leftarrow$(3.1)        &\cite{KL96}\\%
12&(3.$3'$)&$p,p+1$&$p+1=2q$         &o        &$\Leftarrow$(3.3)        &(New) \\%\cite{KL01}
13&(4.2)  &$p,p+1$&$p+1=2q$          &u        &$\Leftarrow$(4.3)        &(New) \\%
14&(4.4)  &$p,p+1$&$p+1=2q$          &u, d     &$\Leftarrow$(4.5)        &(New) \\%
15&(4.6)  &$p,p+1$&$p+1=2q$          &u, d     &$\Leftarrow$(4.2),\,(4.4)&(New) \\%
16&(4.$6'$)&$p,p+1$&$p+1=2q$    &d$\backslash$u&$\Leftarrow$(4.6)        &(New) \\%
\hline
\multicolumn{7}{|c|}{By valency:}\\
\hline
17&(3.1)  &$p,q$  &$p+1=2q$          &u, d     &\quad Analyt.&\cite{KL96}\\%
18&(3.3)  &$p,p+1$&$p+1=2q$          &o        &\quad Analyt.&New        \\%\cite{KL01}
19&(4.3)  &$p,p+1$&$p+1=2q$          &u        &\quad Analyt.&New        \\%
20&(4.$3'$)&$p,p+1$&$p+1=2q$         &u        &$\Leftarrow$(4.3)        &(New) \\%
21&(4.5)  &$p,p+1$&$p+1=2q$          &u, d     &\quad Analyt.&New        \\%
22&(4.7)  &$p,p+1$&$p+1=2q$     &d$\backslash$u&$\Leftarrow$(4.3),\,(4.5)&(New) \\%
23&(3.8)  &$2n$   & $n<27$           &u        &\ Exh.\,search           &\cite{Mc95}\\%9
  &       &       &square-free\,$\rm^c$&       &\quad Analyt.&\cite{KL01}\\%
\hline
\multicolumn{7}{|c|}{Alternating:}\\
\hline
24&(6.1)  &$n$    &$p^2$ or sq. free $\rm^d$&d, sd&\quad Analyt.&New        \\%
25&(6.2)  &$n$    &$p^2$ or sq. free $\rm^d$&u, su&\quad Analyt.&New        \\%
26&(6.3)  &$n$ &$\ p^2$ or sq. free\,$\rm^{d\,e}$&o&\quad Analyt.&New        \\%(6.3a),(6.3b)
27&(6.4)  &$p$    &-                 &u, d     &$\Leftarrow$(6.1),\,(4.1)&(New) \\%(6.2)
\hline
\multicolumn{7}{|c|}{Miscellaneous (non-CI, mixed, of even semi-valency):}\\
\hline
28&(5.2)  &$p,p^2$&-                 &su       &\quad Algebr.&\cite{LP00}\\%{\small non-CI} \cite{KL01}
  &       &       &                  &t        &         &       \\%{\small non-CI}
  &       &       &                  &sd       &         &       \\%{\small non-CI}
29&(5.3)  &$p,p^2$&-                 &su, t, sd&\quad Algebr.&New    \\%
30&(5.4)  &$p^2$  &-                 &su, t, sd&\quad Algebr.&New    \\%{\small mixed}
31&(6.7)  &$p$    &-           &u$^{\rm e}$, sd&$\Leftarrow$(6.1),\,(4.1)&(New) \\%(6.1/2)
\hline
\multicolumn{7}{l}{$\rm^a$ $p$ and $q$ are prime.}\\
\multicolumn{7}{l}{$\rm^b$ Holds also for $n=p^k$ and square-free
$n$ with all prime divisors $p\equiv3$\,(mod\,4).}\\
\multicolumn{7}{l}{\ \ Is conjectured to hold for arbitrary $n$
with all such prime divisors.}\\ %% $p\equiv3$\,(mod\,4).}\\
\multicolumn{7}{l}{$\rm^c$ Is conjectured to hold for arbitrary even orders.}\\
\multicolumn{7}{l}{$\rm^d$ Is conjectured to hold for arbitrary odd orders.}\\
\multicolumn{7}{l}{$\rm^e$ There is a corresponding conjecture for
arbitrary even orders $n,\ 8\nmid n$.}
\end{tabular}

\end{table}}
\end{document}